\documentclass[a4paper,reqno,11pt]{amsart}

\usepackage[utf8]{inputenc}
\usepackage[T1]{fontenc}
\usepackage{lmodern}
\usepackage{amsmath, amsfonts, amssymb, amsthm, mathtools, cancel} 
\usepackage{geometry} 
\usepackage{hyperref} 
\hypersetup{colorlinks = true}
\hypersetup{linktocpage} 
\usepackage[capitalise]{cleveref} 

\newcommand{\interval}[4]{\mathopen{#1}#2
	\mathclose{}\mathpunct{},#3
	\mathclose{#4}}
\newcommand{\intcc}[2]{\interval{[}{#1}{#2}{]}}
\newcommand{\intoc}[2]{\interval{(}{#1}{#2}{]}}
\newcommand{\intco}[2]{\interval{[}{#1}{#2}{)}}
\newcommand{\intoo}[2]{\interval{(}{#1}{#2}{)}}

\newcommand{\abs}[1]{\left\lvert#1\right\rvert}
\newcommand{\norm}[1]{\left\lVert#1\right\rVert}
\newcommand{\scalprod}[2]{\left\langle#1,#2\right\rangle}

\newcommand{\setst}[2]{\left\{#1\mathrel{}\middle|\mathrel{}#2\right\}}

\newcommand{\et}{\quad \text{and} \quad}

\newcommand{\numberset}[1]{\mathbb{#1}}
\newcommand{\N}{\numberset{N}}

\newcommand{\R}{\numberset{R}}
\newcommand{\C}{\numberset{C}}

\newcommand{\coulomb}{\mathcal{C}}
\newcommand{\cutoffsp}{\mathfrak{H}_\mathrm{\Lambda}}
\newcommand{\opspace}{\mathcal{H}_\mathrm{\Lambda}}
\newcommand{\schatten}[1]{\mathfrak{S}_{#1}\left(\mathfrak{H}_\mathrm{\Lambda}\right)}

\newcommand{\diff}{\mathop{}\mathopen{}\mathrm{d}}

\DeclareMathOperator{\tr}{Tr}

\DeclareMathOperator{\supp}{supp}

\theoremstyle{plain}
\newtheorem{thm}{Theorem}[section]
\newtheorem{prop}[thm]{Proposition}
\newtheorem{cor}[thm]{Corollary}
\newtheorem{lem}[thm]{Lemma}

\theoremstyle{definition}

\theoremstyle{remark} 
\newtheorem{rmkx}[thm]{Remark}

\newenvironment{rmk}
  {%
   \pushQED{\qed}\begin{rmkx}}
  {\popQED\end{rmkx}}

\title[Relativistic electrons coupled with Newtonian nuclear dynamics]{Relativistic electrons coupled with Newtonian nuclear dynamics}
\author[Umberto Morellini]{Umberto Morellini}
\address{Ceremade (UMR 7534), Université Paris Dauphine - PSL, Place du Maréchal de Lattre de Tassigny, 75775 Paris Cedex 16, France}
\email{\href{mailto:morellini@ceremade.dauphine.fr}{\nolinkurl{morellini@ceremade.dauphine.fr}}}
\date{25 October 2024}

\begin{document}

\maketitle

\begin{abstract}
In the study of the electronic structure of heavy atoms, relativistic effects cannot be neglected and the Dirac operator naturally appears in place of the Schrödinger operator, raising a number of additional difficulties. The complexity of these systems has been addressed by various approximations.\\
We consider a model consisting of differential equations coupling the time evolution of a finite number of relativistic electrons with the Newtonian dynamics of finitely many nuclei. The electrons are described by the Bogoliubov-Dirac-Fock (BDF) equation of quantum electrodynamics (QED). This system can be seen as a first step in the study of molecular dynamics phenomena, when both relativistic and quantum effects have to be taken into account for the electrons. We address the Cauchy problem for this model and prove global
well-posedness.
\end{abstract}

\tableofcontents

\section{Introduction}\label{sec:introduction}
In quantum chemistry several approximations of the time-independent Schr\"{o}dinger equation have been deeply studied in the past decades, leading to extremely accurate results in the calculation of the ground state of molecules containing several hundred electrons. Nevertheless, it is quite evident that these models are insufficient to study dynamical phenomena such as chemical reactions which rule the behaviour of matter. For this reason, quantum chemists have been encouraged to study time-dependent models.\\
Whenever relativistic effects have to be taken into account, the Dirac operator naturally appears in place of the Schr\"{o}dinger operator, raising a number of additional difficulties which will be profoundly discussed later on. As before, approximations become necessary in order to deal with the complexity of the system. The main purpose of this paper is to study the global well-posedness of the Cauchy problem associated with a Bogoliubov-Dirac-Fock dynamics of a finite number of relativistic electrons coupled with classical dynamics of finitely many nuclei. This model can be considered as a step towards a better comprehension of dynamical phenomena when both relativistic and quantum effects have to be taken into account. This work has been inspired by paper \cite{CanLeB-1999-MMMAS} where the authors consider a non-adiabatic approximation of the Schr\"{o}dinger equation such that the electrons are ruled by time-dependent Hartree-Fock equations and the nuclei-electrons interaction is of Hellmann-Feynman type. For an introduction about this class of potentials (or forces), we refer the reader to \cite[Section~II]{BorNetSch-1996-JCP} and \cite[Section~2.3]{CanLeB-1999-MMMAS}.\\
In this section, we introduce the free Dirac operator to describe the relativistic dynamics of electrons and discuss the mathematical and physical difficulties carried by this operator. This leads us to talk about the Bogoliubov-Dirac-Fock (BDF) model, which is a mean-field approximation of quantum electrodynamics (QED) allowing us to overcome all these problems. In particular, this model allows us to consider electrons interacting with the quantum vacuum in addition to the electron-electron and nucleus-electron interactions. Indeed, in QED the vacuum is not an empty space, but rather a fluctuating medium behaving as a non-linear polarisable material. Later we discuss the dynamics of the nuclei, by justifying the classical approximation of their nature that we take into account.\\
In \cref{sec:the-general-setting}, we introduce our model describing a finite number of relativistic electrons obeying the time-dependent BDF equation and interacting with the Newtonian dynamics of a finite number of nuclei whose centres of mass are considered to be classical and interacting themselves with the relativistic electrons. We state our main result in the most general form, but then the proof presented in the following sections is restricted to the case with only two nuclei for the sake of simplicity. It is important to notice that our argument applies to any system consisting of a finite number of nuclei and electrons without any remarkable change. In the one nucleus case, the nucleus-nucleus interaction term can be neglected and computations are slightly simpler.\\
In \cref{sec:local-existence} and \cref{sec:uniqueness}, we prove a local existence and uniqueness result for the solution of the two-nuclei model by the Schauder fixed-point theorem and a standard application of Gr\"{o}nwall's lemma. Here it is important to define a suitable function space to apply a fixed-point argument.\\
In \cref{sec:global-in-time-existence}, we extend the results of the previous sections by proving the global-in-time existence of the maximal solution by means of the conservation of energy of the system. Here the pattern of the proof follows a standard argument used in \cite{HaiLewSpa-2005-LMP}.\\
It is important that we assume that no collisions happen between nuclei throughout this work since the conservation of energy does not prevent them a priori. Indeed, a proof of this result would require some additional estimates on the vacuum-nuclei Coulomb interaction term. It would be interesting to study this problem in more detail in the future. Moreover, notice that due to the assumptions our model cannot cover some physical situations. For instance, it is well known that the motion of charged nuclei creates magnetic fields which are not taken into account by this model. In this case, the situation becomes much more complicated and different models giving a better approximation should be addressed. In particular, more sophisticated regularisation techniques need to be considered in order to deal with ultraviolet divergence by preserving the gauge invariance as explained in \cite{GraHaiLewSer-2013-ARMA,GraLewSer-2018-JMPA}. This will be the object of future investigations.

\subsection{The Bogoliubov-Dirac-Fock model}
The dynamics of relativistic quantum spin-$1/2$ particles is described by the free Dirac operator,
\[
D^0=-i\boldsymbol{\alpha}\cdot\nabla+\beta=-i\sum_{k=1}^3
\alpha_k\partial_k+\beta,
\]
acting on $L^2\left(\R^3,\C^4\right)$, introduced by Dirac in $1928$ \cite{Dir-1928-PRSL}, where $\boldsymbol{\alpha}=\left(\alpha_1,\alpha_2,\alpha_3\right)$ and
\[
\beta=\begin{bmatrix}
I_2 & 0\\
0 & I_2
\end{bmatrix},\quad
\alpha_k=\begin{bmatrix}
0 & \sigma_k\\
\sigma_k & 0
\end{bmatrix},
\]
with
\[
\sigma_1=\begin{bmatrix}
0 & 1\\
1 & 0
\end{bmatrix},\quad
\sigma_2=\begin{bmatrix}
0 & -i\\
i & 0
\end{bmatrix},\quad
\sigma_3=\begin{bmatrix}
1 & 0\\
0 & -1
\end{bmatrix}.
\]
We point out that we choose a system of units such that $\hbar=c=1$ and the mass of the electron is equal to $1$ as well. We refer to \cite{Tha-1992-book} for a detailed reading about this operator and its properties, which have been and are currently deeply studied. For our purpose, it is enough to know that its spectrum is not bounded from below. Indeed,
\begin{equation}\label{eq:free-dirac-spectrum}
\sigma\left(D^0\right)=\intoc{-\infty}{-1}\cup\intco{1}{+\infty}.
\end{equation}
This leads to important mathematical, physical and numerical difficulties. In particular, \eqref{eq:free-dirac-spectrum} led Dirac to postulate that the vacuum is filled with infinitely many particles occupying the negative energy states (Dirac sea) \cite{Dir-1930-PRSL, Dir-1934-MPCPS, Dir-1934-SR}. The distribution of negative energy electrons is postulated to be unobservable because of its uniformity. As a consequence, a real free electron cannot be in a negative energy state due to the Pauli principle. Moreover, Dirac conjectured the existence of "holes" in the Dirac sea interpreted as positrons, having a positive charge and a positive energy. Dirac also predicted the phenomenon of vacuum polarisation: in the case of an external electric field, the negative energy electrons align their charges with the field direction and the vacuum obtains a non-constant density of charge. Actually the polarised vacuum modifies the electrostatic field and the virtual electrons react to the corrected field.\\
The vacuum polarisation is almost negligible in the calculation of Lamb shift for the ordinary hydrogen atom, but it is crucial for high-$Z$ atoms \cite{MohPluSof-1998-PR}. Indeed, for heavy atoms it is necessary to take relativistic effects into account. Nevertheless, we are not able to define the equivalent of the $N$-body non-relativistic Hartree-Fock theory by replacing the Schr\"{o}dinger operator with the Dirac operator, because of the negative part of its spectrum. The correct theory is quantum electrodynamics which is the quantum field theory of electrons and positrons (the electron-positron field) and photons (the electromagnetic field). Its predictive success about the interaction between matter and light relies on pertubation theory. However, its description in terms of perturbative calculations restricts the range of physical situations where it can be applied. A mathematically consistent non-perturbative formulation of QED is still unknown and the argumentation above makes it an interesting open problem. To overcome the difficulty, several mean-field models have been proposed and studied in the past decades in order to get a first approximation. In particular, effective models deduced from non-relativistic theories (such as the Dirac-Fock model) suffer from inconsistencies: for instance, a ground state never minimises the physical energy which is always unbounded from below. We refer the reader to \cite{EstLewSer-2008-BAMS} for a complete review.\\
The phenomenon of vacuum polarisation described above suggests that mathematically speaking we need to find a self-consistent equation for the polarised vacuum and to solve it by a standard fixed-point argument. In physics, self-consistent equations are usually derived as Euler-Lagrange equations of an energy functional. For instance, this is the case in the non-relativistic Hartree-Fock model. Instead, the same argument does not apply to the Dirac-Fock model, precisely because of the unboundedness from below of the spectrum of the Dirac operator. For this reason, we need to study the Bogoliubov-Dirac-Fock energy functional, first introduced in $1989$ \cite{ChaIra-1989-JPB}. They derive their mean-field model from the QED Hamiltonian by neglecting the contribution of photons (the electromagnetic field is not quantised) and normal-ordering it with respect to the free projector,
\[
P^0=\chi_{\left(-\infty,0\right)}\left(D^0\right),
\]
as a fixed reference. They do not neglect the vacuum polarisation terms by pointing out their importance for the consistency of the relativistic mean-field theory. Here, $P^0$ formally corresponds to the one-particle density matrix of an infinite Slater determinant. This model describes the behaviour of a finite number of spin-$1/2$ particles coupled to the Dirac sea which can become polarised. Indeed, the introduction of the vacuum in the model is the solution to deal with the negative energies of the Dirac operator and obtain a bounded-below energy functional, whose minimisers are well defined ground states. Finally we get a model which is more accurate (taking into account vacuum polarisation effects), more physically sensible and whose ground states have a proper definition. On the other hand, one has to deal with infinitely many particles (the real ones and the virtual ones of the Dirac sea) instead of the finite number of electrons. Moreover, it is not easy to give a meaning to the quantities involved in our model because several divergent quantities appear in QED even after normal ordering. We refer the reader to \cite{HaiLewSerSol-2007-PRA} for a brief and precise review about the argument.\\
Mathematically speaking, for any $Q=P-P^0\in\mathfrak{S}_2\left(L^2\left(\R^3,\C^4\right)\right)$ ($P$ being an orthogonal projector), the BDF energy can be formally written as follows:
\begin{equation}\label{eq:BDF-energy}
E\left(Q\right)=\tr\left(D^0Q\right)-\alpha\int\rho_Q\varphi+\frac{\alpha}{2}\int\int\frac{\rho_Q\left(x\right)\rho_Q\left(y\right)}{\abs{x-y}}\diff x\diff y-\frac{\alpha}{2}\int\int\frac{\abs{Q\left(x,y\right)}^2}{\abs{x-y}}\diff x\diff y,
\end{equation}
where $\rho_Q\left(x\right)=\tr_{\C^4}Q\left(x,x\right)$ is the charge density associated to $Q$ and $\mathfrak{S}_2\left(L^2\left(\R^3,\C^4\right)\right)$ denotes the class of Hilbert-Schmidt operators on $L^2\left(\R^3,\C^4\right)$. Here, $\varphi$ is the external electrostatic potential created by the (non-negative) charge distribution $f$, namely $\varphi=f*|\cdot|^{-1}$. By formal computation, we can show that the Euler-Lagrange equation of this functional is
\begin{equation}\label{eq:stationnary-BDF}
\left[P,D_Q\right]=0,
\end{equation}
where
\[
D_Q=D^0-\alpha\varphi+\alpha\rho_Q*\frac{1}{|\cdot|}-\alpha\frac{Q\left(x,y\right)}{\abs{x-y}}
\]
is the so called mean-field operator taking into account not only external fields, but also the self-consistent potential created by the state $P\left(t\right)$ itself. In the time-dependent case, \cref{eq:stationnary-BDF} modifies into
\[
i\frac{\diff}{\diff t}P\left(t\right)=\left[D_{Q\left(t\right)},P\left(t\right)\right].
\]
We have to deal with several difficulties in order to well define the BDF energy in \eqref{eq:BDF-energy}. First of all, we have to deal with divergences in the high-frequencies domain to be able to define, for instance, the density $\rho_Q$. In this case, one possible solution is to introduce an ultraviolet cut-off $\mathrm{\Lambda}>0$, by replacing the ambient space $L^2\left(\R^3,\C^4\right)$ with
\[
\cutoffsp=\setst{f\in L^2\left(\R^3,\C^4\right)}{\supp\widehat{f}\subseteq B\left(0,\mathrm{\Lambda}\right)}.
\]
Secondly, we need to find a way to define $\tr\left(D^0Q\right)$ and for this reason the definition of trace functional has to be extended to a bigger class of compact operators. In \cite{HaiLewSer-2005-CMP}, the authors define an operator $A\in\mathfrak{S}_2\left(\cutoffsp\right)$ to be $P^0$-trace class if $A^{++}=\left(1-P^0\right)A\left(1-P^0\right)\in\mathfrak{S}_1\left(\cutoffsp\right)$ and $A^{--}=P^0AP^0\in\mathfrak{S}_1\left(\cutoffsp\right)$ (the set of all trace class operators). The associated $P^0$-trace of $A$ is then given by
\[
\tr_{P^0}A=\tr A^{++}+\tr A^{--}.
\]
The set of all $P^0$-trace class operators is denoted by $\mathfrak{S}^{P^0}_1\left(\cutoffsp\right)$. Notice that we do not assume the initial condition to be trace class, since this property is not conserved along the solution, namely $Q\left(t\right)\notin\mathfrak{S}_1\left(\cutoffsp\right)$ as soon as $t>0$, even if the initial datum $Q\left(0\right)$ is assumed to be $\mathfrak{S}_1\left(\cutoffsp\right)$.\\
It was proved in \cite[Lemma~2]{HaiLewSer-2005-CMP} that any difference of two projectors satisfying $Q=P-P^0\in\mathfrak{S}_2\left(\cutoffsp\right)$ (Shale-Stinespring criterion) is automatically in $\mathfrak{S}^{P^0}_1\left(\cutoffsp\right)$ and moreover $\tr_{P^0}\left(Q\right)$ is always an integer. Thus the $P^0$-trace is a suitable tool for describing charge sectors, which is crucial when we want to work with physical systems containing also a finite number $q$ of real electrons in addition to the infinite virtual ones. Indeed, in this case we need to minimise the well defined BDF energy in a fixed charge sector, i.e. under the constraint $\tr_{P^0}\left(Q\right)=q$. Later in this paper, we will see that this quantity is conserved all along the BDF dynamics. This allows us to study the interaction between a finite number of relativistic electrons and a finite number of classical nuclei.\\
We refer the reader to \cite{HaiLewSer-2005-CMP, HaiLewSer-2005-JPA, HaiLewSol-2007-CPAM, GraLewSer-2009-CMP, HaiLewSer-2009-ARMA} for all the rigorous details and main results about the BDF model.

\subsection{The classical dynamics for nuclei}
The Born-Oppenheimer approximation is one of the most known and useful mathematical approximations to deal with molecular dynamics in quantum chemistry. It allows simplifying the Schr\"{o}dinger equation by treating separately the wavefunctions of atomic nuclei and electrons. However there are physical situations where this approximation does not work, like for instance in presence of a time-dependent electric field. We refer the reader to \cite[Paragraphs~2.2 and 2.3]{CanLeB-1999-MMMAS} for a more detailed discussion about the adiabatic approximation, which is a generalisation of the Born-Oppenheimer approximation, and its range of applications.\\
In order to deal with this kind of situations, we need to introduce another model to study the dynamics of the nuclei. Since it is well known in the literature that point-like sources create divergences in classical electrodynamics as well as in QED \cite{LanAbrHal-1956-INC}, we take nuclear charge distributions centred in the centres of mass of the corresponding nuclei and suppose them to be classical. The physical explanation of this approximation relies on the fact that atomic nuclei are much heavier than electrons. Indeed, the ratio is approximately $1836$ for the hydrogen nucleus and is greater than $10^4$ for most of the other atoms encountered in chemistry. Hence, the quantum nature of the nuclei can be neglected in most applications. This approximation is almost always valid in chemistry and allows us to study the time evolution of atomic nuclei from a classical point of view. Moreover, the dynamics of the classical nuclei can be assumed to be non-relativistic mainly because these particles are really massive and therefore heavy. For this reason, they are usually slow compared to the speed of light and their relativistic corrections can be neglected.\\
Let us consider $M$ nuclei and denote by $m_k$ and $z_k$ respectively the mass and the charge of the $k$-th nucleus. The state of the system is then classically described at time $t$ as a point in the phase space given by the coordinates of the centres of mass,
\[
\left(x_k\left(t\right),\frac{\diff x_k\left(t\right)}{\diff t}\right)_{1\leq k\leq M}\in\R^{6M}.
\]
Moreover we assume that the $k$-th nucleus carries a density of charge $f_k$ which creates the corresponding electrostatic potential $\varphi_k\left(t\right)=z_kf_k\left(\abs{\cdot-x_k\left(t\right)}\right)*\frac{1}{\abs{\cdot}}$ at time $t$. We do not need to fix the sign of the charge density from a mathematical point of view, but we can think about a non-negative function in order to represent for instance a system of nuclei in a molecule. We assume that $f_k\in L^1\left(\R^3,\R\right)\cap\coulomb$, $1\leq k\leq M$, where
\[
\coulomb=\setst{\rho\in\mathcal{S}'\left(\R^3\right)}{\widehat{\rho}\in L^1_{\mathrm{loc}}\left(\R^3\right),\mathcal{D}\left(\rho,\rho\right)<\infty}
\]
is the Coulomb space endowed with the following scalar product:
\[
\mathcal{D}\left(\rho_1,\rho_2\right)=\int_{\R^3}\frac{\overline{\widehat{\rho}_1\left(k\right)}\widehat{\rho}_2\left(k\right)}{\abs{k}^2}\diff k.
\]
We immediately remark that $\coulomb=\dot{H}^{-1}\left(\R^3\right)$. We need the $L^1$-integrability in order to suppose the functions to be normalised (which is physically reasonable), namely
\[
\int_{\R^3}f_k\left(x\right)\diff x=1,\quad 1\leq k\leq M.
\]
\begin{rmk}
Actually, later in this paper, the charge distributions $f_k$ are required to be more regular. Indeed, in \cref{stat:main-result} the distribution of charge is needed to belong additionally to $L^2\left(\R^3\right)$. Mathematically speaking this is a technical constraint, but from a computational point of view this is not excessively restrictive since quantum chemists are used to choose nuclear charge distributions which are smooth as explained in \cite{VisDya-1997-ADNDT}.
\end{rmk}

\section{The general setting and main result}\label{sec:the-general-setting}
We are now ready to introduce the complete model which is studied in this work. We work on relativistic electrons interacting with a finite number $M$ of classical nuclei of charge $z_k$ and mass $m_k$, $k=1,\ldots M$. Let us call $f_k$ the corresponding charge distributions with the regularity properties described in the previous section. The Bogoliubov-Dirac-Fock dynamics of any finite number of relativistic electrons is described by the following equation:
\begin{equation}\label{eq:time-dependent-BDF}
i\frac{\diff}{\diff t}P\left(t\right) = \left[D_{Q,x_1,\dots,x_M},P\left(t\right)\right],
\end{equation}
where the mean-field operator at time $t$,
\[
D_{Q,x_1,\dots,x_M}=\mathcal{P}_\mathrm{\Lambda}\left(D^0-\alpha\sum_{k=1}^M z_k f_k\left(\abs{\cdot-x_k}\right)*\frac{1}{|\cdot|}+\alpha\rho_Q*\frac{1}{|\cdot|}-\alpha\frac{Q\left(x,y\right)}{\abs{x-y}}\right)\mathcal{P}_\mathrm{\Lambda}
\]
is a well defined self-adjoint bounded operator on $\cutoffsp$. Notice that without the orthogonal projection $\mathcal{P}_\mathrm{\Lambda}$ of $L^2\left(\R^3\right)$ onto $\cutoffsp$ the operator $D_{Q,x_1,\dots,x_M}$ would not preserve the Hilbert space $\cutoffsp$. We remark that the nuclear coordinates $x_k \left(t\right)$ are here parameters for the mean-field operator. In several context, we can just get rid of them, by neglecting the motion of the nuclei: therefore the nuclei do not move and we only have to deal with \cref{eq:time-dependent-BDF}. Since we are interested in describing situations where the dynamics of the nuclei is important, such as chemical reactions, we couple \eqref{eq:time-dependent-BDF} with
\begin{equation}\label{eq:newton}
m_k\frac{\diff^2}{\diff t^2}x_k\left(t\right) = -\nabla_{x_k} W_Q\left(t,x_1,\dots,x_M\right),
\end{equation}
where
\begin{multline*}
    W_Q\left(t,x_1,\dots,x_M\right)=\alpha \mathcal{D}\left(\rho_Q,\sum_{k=1}^M z_k f_k\left(\abs{\cdot-x_k}\right)\right)\\
    -\alpha\sum_{1\leq i<j\leq M}\int \int \frac{z_i f_i\left(\abs{x-x_i}\right) z_j f_j\left(\abs{y-x_j}\right)}{\abs{x-y}}\diff x\diff y
\end{multline*}
takes into account the contribution given by the scalar product defined in the Coulomb space $\coulomb$ describing the Coulomb interaction between electrons and nuclei and the nucleus-nucleus interactions. \cref{eq:newton} means that in the approximation discussed above the centre of mass of the nucleus moves obeying the Newtonian dynamics in the electrostatic potential created by electrons and nuclei.
We now introduce our main Hilbert space,
\[
\opspace=\setst{Q\in\schatten{2}}{\rho_Q\in\coulomb},
\]
endowed with the following norm:
\[
\norm{Q}=\left(\norm{Q}^2_{\schatten{2}}+\norm{\rho_Q}^2_\coulomb\right)^\frac{1}{2}.
\]
Thus, we finally consider the following:
\begin{equation}\label{eq:M-nuclei}
\begin{cases}
i\frac{\diff}{\diff t}P\left(t\right) = [D_{Q,x_1,\dots,x_M},P\left(t\right)],\\[8pt]
m_k\frac{\diff^2}{\diff t^2}x_k\left(t\right) = -\nabla_{x_k}W_Q\left(t,x_1,\dots,x_M\right),\quad k=1,\ldots M\\[8pt]
P\left(0\right) = P_I,\;P\left(t\right)^2 = P\left(t\right),\;Q\left(t\right) = P\left(t\right) - P\left(0\right) \in\opspace,\\[8pt]
x_k\left(0\right) = x^0_k\in\R^3,\;\frac{\diff x_k}{\diff t}\left(0\right) = v^0_k\in\R^3,\quad k=1,\ldots M
\end{cases}
\end{equation}
By combining the well defined BDF energy,
\[
E^{\mathrm{BDF}}\left(Q\right)=\tr_{P^0}\left(D^0Q\right)-\alpha \mathcal{D}\left(\rho_Q,\sum_{k=1}^M z_k f_k\right)+\frac{\alpha}{2}\mathcal{D}\left(\rho_Q,\rho_Q\right)-\frac{\alpha}{2}\int\int_{\R^6}\frac{\abs{Q\left(x,y\right)}^2}{\abs{x-y}}\diff x\diff y,
\]
(acting on $Q=P-P^0\in\opspace$, $P$ being an orthogonal projector) with the well known Newtonian energy, we can define the energy of the system described by \eqref{eq:M-nuclei} as follows:
\begin{multline*}
E^{\left(M\right)}\left(Q\left(t\right),x_1\left(t\right),\ldots,x_M\left(t\right)\right)\\
=\tr_{P^0}\left(D^0Q\left(t\right)\right)-\alpha \mathcal{D}\left(\rho_{Q\left(t\right)},\sum_{k=1}^M z_k f_k\left(\abs{\cdot-x_k\left(t\right)}\right)\right)+\frac{\alpha}{2}\mathcal{D}\left(\rho_{Q\left(t\right)},\rho_{Q\left(t\right)}\right)\\
-\frac{\alpha}{2}\int\int\frac{\abs{Q\left(x,y\right)}^2}{\abs{x-y}}\diff x\diff y+\frac{1}{2}\sum_{k=1}^M m_k \abs{\dot{x}_k\left(t\right)}^2\\
+\alpha\sum_{1\leq i<j\leq M}\int \int \frac{z_i f_i\left(\abs{x-x_i\left(t\right)}\right) z_j f_j\left(\abs{y-x_j\left(t\right)}\right)}{\abs{x-y}}\diff x\diff y.
\end{multline*}
In \cite{CacDeSNoj-2019-MA}, the authors introduce the relativistic dynamics in the model of \cite{CanLeB-1999-MMMAS}, by replacing the Schr\"odinger operator by the Dirac operator. In this way, they describe the time evolution of one relativistic electron interacting with finitely many nuclei, without considering the effects of the Dirac sea interacting with the nuclear and electronic external electric fields. They prove a local existence result under a constraint on the total charge of the nuclei. Our main purpose here consists of extending both the result proved in \cite{HaiLewSpa-2005-LMP} by including the nuclear dynamics into the time evolution of the system and the result proved in \cite{CanLeB-1999-MMMAS} by taking into account relativistic quantum effects such as vacuum polarisation. Therefore we want to show the global well-posedness of Cauchy problem~\eqref{eq:M-nuclei}, namely that the system has a unique global solution in a suitable function space which is discussed later on. We can now state our main result whose proof is detailed in the sequel of this paper.
\begin{thm}\label{stat:main-result}
Let $0\leq\alpha<4/\pi$ and $f_k\in L^1\left(\R^3\right)\cap L^2\left(\R^3\right)\cap\coulomb$, $k=1,\ldots M$. Let $P_I$ be an orthogonal projector such that $Q_I=P_I-P^0\in\opspace$ and $x^0_k,v^0_k\in\R^3$, $k=1,\ldots M$. Then there exists a unique global solution
\[
\left(Q,x_1,\ldots,x_M\right)\in C^1\left(\intco{0}{+\infty},\opspace\right)\times \left(C^2\left(\intco{0}{+\infty},\R^3\right)\right)^M
\]
of system~\eqref{eq:M-nuclei}. Moreover, $Q\left(t\right)=P\left(t\right)-P^0$ is $P^0$-trace class and
\[
\tr_{P^0}\left(Q\left(t\right)\right)=\tr_{P^0}\left(Q_I\right)
\]
and
\[
E^{\left(M\right)}\left(Q\left(t\right),x_1\left(t\right),\ldots,x_M\left(t\right)\right)=E^{\left(M\right)}\left(Q_I,x^0_1,\ldots,x^0_M\right),
\]
for all $t\in\intco{0}{+\infty}$.
\end{thm}

\section{Local existence}\label{sec:local-existence}
As said in the introduction, for the sake of simplicity we prove the main result of this work in the case with only two nuclei. For an arbitrary number of nuclei, the proof works exactly as in the presence of just two nuclei. Indeed, concerning the BDF equation of the system, the only different term is the nuclei-electrons Coulomb scalar product which contains as many terms as nuclei ($\mathcal{D}\left(\rho_Q,z_k f_k\right)$, $k=1,\ldots,M$), each of which can be estimated in the same way. Concerning the Newton equation, the potential $W_Q$ contains the electrons-nuclei interactions as well as the nuclei-nuclei interactions, which can be estimated in the same way independently of the number of nuclei. In the one-nucleus case, the nuclei-nuclei interactions disappear from the Newton equation for obvious reasons and the calculations become slightly simpler, as already mentionned in the introduction. For all of these reason, the statement can be proved without loss of generality in the presence of just two nuclei of mass $m_1,m_2$, charge $z_1,z_2$ and charge distribution $f_1,f_2$ satisfying the hypotheses discussed above. We start by proving a local existence result for solutions of the evolution problem \eqref{eq:M-nuclei} written in terms of $Q\left(t\right)=P\left(t\right)-P^0$, namely
\begin{equation}\label{eq:2-nuclei}
\begin{cases}
i\frac{\diff}{\diff t}Q\left(t\right)=F\left(Q\right),\\[8pt]
m_k\frac{\diff^2}{\diff t^2}x_k\left(t\right) = -\nabla_{x_k}W_Q\left(t,x_1,x_2\right),\quad k=1,2\\[8pt]
Q\left(0\right)=Q_I\coloneqq P_I-P^0\in\opspace,\\[8pt]
x_k\left(0\right) = x_k^0\in\R^3,\;\frac{\diff x_k}{\diff t}\left(0\right) = v_k^0\in\R^3,\quad k=1,2
\end{cases}
\end{equation}
where
\[
F(Q)=\left[D_{Q,x_1,x_2},Q\right]+\left[V_{Q,x_1,x_2},P^0\right],
\]
acting on $\opspace$, with
\[
V_{Q,x_1,x_2}=\mathcal{P}_\mathrm{\Lambda}\left(\alpha\left(\rho_Q-z_1 f_1\left(\abs{\cdot-x_1}\right)-z_2 f_2\left(\abs{\cdot-x_2}\right)\right)*\frac{1}{|\cdot|}-\alpha\frac{Q\left(x,y\right)}{\abs{x-y}}\right)\mathcal{P}_\mathrm{\Lambda}.
\]
In order to do that, we want to use a Schauder fixed-point argument and it is then important to define suitable function spaces where correctly applying our machinery.\par
We denote by $B_{\opspace}\left(0,\norm{Q_I}\right)$ the ball of radius $\norm{Q_I}$ centered at the origin in the operator space $\opspace$, equipped with the corresponding topology. Then, $F$ is continuous on $B_{\opspace}\left(0,\norm{Q_I}\right)$ with respect to the $\opspace$ topology (we refer the reader to Lemma~\ref{stat:BDF-existence} below for a proof of its continuity). Indeed, because of the presence of a quadratic term, $F$ needs to be considered on a bounded space in order to be continuous. Let us denote $C_F$ the corresponding continuity constant of $F$, which will depend on the norm of the initial datum $\norm{Q_I}$. We can then fix some arbitrary time $\tau>0$ such that
\begin{gather}
\label{eq:inequality-1}
\frac{1}{1-\tau C_F}\leq 1,\\
\label{eq:inequality-2}
\alpha\frac{\tau}{m_k}C_f\left(\norm{Q_I}+C_f\right)\leq 1\quad k=1,2,
\end{gather}
where
\[
C_f =\max_{k=1,2}\{\sup_{\R^3}\norm{z_k \nabla f_k\left(\abs{\cdot-x}\right)}_\coulomb,\sup_{\R^3}\norm{z_k f_k\left(\abs{\cdot-x}\right)}_\coulomb\}.
\]
Notice that $\tau$ will then depend on the norm of the initial datum and not on the initial datum itself. Now, we define the sets:
\[
\mathcal{B}^e_\tau=\setst{Q\in C^1\left(\intcc{0}{\tau},\opspace\right)}{\norm{Q}_{C^0\left(\intcc{0}{\tau},\opspace\right)}\leq \norm{Q_I}},
\]
\begin{multline*}
    \mathcal{B}^{n,k}_\tau=\setst{x_k\in C^1\left(\intcc{0}{\tau},\R^3\right)}{x_k\left(0\right)=x_k^0,\,\frac{\diff x_k}{\diff t}\left(0\right)=v_k^0,\,\norm{\frac{\diff x_k}{\diff t}}_{C^0(\intcc{0}{\tau},\R^3)}\leq \abs{v_k^0}+1},\\
    k=1,2,
\end{multline*}
endowed with the topology of $C^0\left(\intcc{0}{\tau},\opspace\right)$ and the topology of $C^1\left(\intcc{0}{\tau},\R^3\right)$ respectively. In the sequel, we shall need to consider also the sets $\mathcal{B}^{n,k}_\tau\cap C^2\left(\intcc{0}{\tau},\R^3\right)$, $k=1,2$, equipped with the topology of $C^2\left(\intcc{0}{\tau},\R^3\right)$.

\subsection{Local existence for the BDF equation}
We now decouple the differential equations in \eqref{eq:2-nuclei} and focus on the Cauchy problem of the first one of them,
\begin{equation}\label{eq:BDF-subsystem}
	\begin{cases}
	i\frac{\diff}{\diff t}Q\left(t\right)=\left[D_{Q,x_1,x_2},Q\right]+\left[V_{Q,x_1,x_2},P^0\right],\\[8pt]
	Q\left(0\right)=Q_I=P_I-P^0\in\opspace.
	\end{cases}
\end{equation}
Before stating and proving the local existence for solutions of \eqref{eq:BDF-subsystem}, let us state the following auxiliary lemma.
\begin{lem}\label{stat:f-convergence}
Let $g\in\coulomb$ and $\left(x_n\right)_{n\in\N}\subseteq\R^3$ be a sequence converging to $x\in\R^3$. Then $\lim_{n\to\infty}g\left(\abs{\cdot-x_n}\right)=g\left(\abs{\cdot-x}\right)$ in $\coulomb$.
\end{lem}
\proof
By the following calculation, one easily has
\begin{align*}
\norm{g\left(\abs{\cdot-x_n}\right)-g\left(\abs{\cdot-x}\right)}_\coulomb^2 &=\int\frac{\abs{\widehat{g\left(|\cdot|\right)}\left(k\right)e^{i\scalprod{x_n}{k}}-\widehat{g\left(|\cdot|\right)}\left(k\right)e^{i\scalprod{x}{k}}}^2}{\abs{k}^2}\diff k\\
    &=\int\frac{\abs{\widehat{g\left(|\cdot|\right)}\left(k\right)}^2}{\abs{k}^2}\abs{e^{i\scalprod{x_n}{k}}-e^{i\scalprod{x}{k}}}^2\diff k.
\end{align*}
Let us now denote
\[
\xi_n\left(k\right)=\frac{\abs{\widehat{g\left(|\cdot|\right)}\left(k\right)}^2}{\abs{k}^2}\abs{e^{i\scalprod{x_n}{k}}-e^{i\scalprod{x}{k}}}^2.
\]
The function $\xi_n$ can be easily dominated by $4 \abs{\widehat{g\left(|\cdot|\right)}\left(k\right)}^2/\abs{k}^2$ which is integrable by hypothesis, uniformly with respect to $n$. Moreover, $\xi_n$ converges to zero pointwise, as $n\to\infty$, by hypothesis on the sequence $\left(x_n\right)_{n\in\N}\subseteq\R^3$. The dominated convergence theorem allows us to conclude the proof.
\endproof
We can now state the following lemma:
\begin{lem}\label{stat:BDF-existence}
Let $x_1\in\mathcal{B}^{n,1}_\tau$ and $x_2\in\mathcal{B}^{n,2}_\tau$. Then system~\eqref{eq:BDF-subsystem} has a unique (local) solution $Q$ in $C^1\left(\intcc{0}{\tau},\opspace\right)$ which is also in $\mathcal{B}^e_\tau$. Moreover, the map
\[
\mathcal{F}^{\mathrm{BDF}}:\left(x_1,x_2\right)\in\mathcal{B}^{n,1}_\tau\times\mathcal{B}^{n,2}_\tau\longmapsto Q\in\mathcal{B}^e_\tau
\]
is bounded and continuous.
\end{lem}	
\proof
By the Cauchy-Lipschitz theorem we need to show that $F$ is locally Lipschitz for the $\opspace$ topology in order to have existence and uniqueness of a local solution of \eqref{eq:BDF-subsystem}. Let $Q_1,Q_2\in\opspace$. Hence,
\[
F\left(Q_1\right)-F\left(Q_2\right)
=\left[V'_{Q_1-Q_2},P^0\right]+\left[D_{x_1,x_2},Q_1-Q_2\right]+\left[V'_{Q_1}-V'_{Q_2},Q_1\right]+\left[V'_{Q_2},Q_1-Q_2\right],
\]
where we have introduced the following notations:
\begin{gather}
\notag
V'_Q=\mathcal{P}_\mathrm{\Lambda}\left(\alpha\rho_Q*\frac{1}{|\cdot|}-\alpha\frac{Q\left(x,y\right)}{\abs{x-y}}\right)\mathcal{P}_\mathrm{\Lambda},\\
\notag
D_{x_1,x_2}=\mathcal{P}_\mathrm{\Lambda}\left(D^0-\alpha z_1 f_1\left(\abs{\cdot-x_1}\right)*\frac{1}{|\cdot|}-\alpha z_2 f_2\left(\abs{\cdot-x_2}\right)*\frac{1}{|\cdot|}\right)\mathcal{P}_\mathrm{\Lambda}.
\end{gather}
So we get four different terms which can all be estimated exactly as in the proof of \cite[Lemma~3.1]{HaiLewSpa-2005-LMP}. In particular, there exist constants $C_1,C_2$ and $C_3$ (dependent on the cut-off parameter $\mathrm{\Lambda}$ and on the Coulomb norms of the external charge distributions $f_1,f_2$) such that
\[
\norm{\left[D_{x_1,x_2},Q\right]}\leq C_1 \norm{Q},\quad\norm{\left[V'_Q,Q'\right]}\leq C_2\norm{Q}\norm{Q'}\et\norm{\left[V'_Q,P^0\right]}\leq C_3 \norm{Q}.
\]
Notice that this proves also the continuity of $F$ on $B_{\opspace}\left(0,C^e\norm{Q_I}\right)$ for the same topology.\\
Now we need to prove that the application $\mathcal{F}^{\mathrm{BDF}}:\left(x_1,x_2\right)\in\mathcal{B}^{n,1}_\tau\times\mathcal{B}^{n,2}_\tau\mapsto Q\in\mathcal{B}^e_\tau$ is bounded and continuous. First of all, the map $\mathcal{F}^{\mathrm{BDF}}$ is well-defined. Indeed,
\[
\norm{Q\left(t\right)}=\norm{Q_I}+\int_0^t\norm{\frac{\diff Q}{\diff s}}\diff s
\leq\norm{Q_I}+\tau\norm{F\left(Q\right)}
\leq\norm{Q_I}+\tau C_F \norm{Q\left(t\right)}
\]
which easily implies $\norm{Q}_{C\left(\intcc{0}{\tau},\opspace\right)}\leq \norm{Q_I}$ by \eqref{eq:inequality-1}.
Of course the application $\mathcal{F}^{\mathrm{BDF}}$ is bounded since $\mathcal{B}^e_\tau$ is bounded.\\
Let us now prove the continuity. Let $\left(x^n_1,x^n_2\right)_{n\in\N} \subseteq\mathcal{B}^{n,1}_\tau\times\mathcal{B}^{n,2}_\tau$ a sequence converging to $\left(x_1,x_2\right)\in\mathcal{B}^{n,1}_\tau\times\mathcal{B}^{n,2}_\tau$. Let us denote $Q^n=\mathcal{F}^{\mathrm{BDF}}\left(x^n_1,x^n_2\right)$ and $Q=\mathcal{F}^{\mathrm{BDF}}\left(x_1,x_2\right)$. We want to prove that $Q^n$ converges to $Q$ in $\mathcal{B}^e_\tau$. Thus,
\[
\norm{Q^n-Q}\leq\int_0^t\norm{\frac{\diff}{\diff\tau}\left(Q^n-Q\right)}\diff\tau\leq\int_0^t\norm{F\left(Q^n\right)-F\left(Q\right)}\diff\tau.
\]
Since $Q^n$ and $Q$ correspond respectively to $x^n_1,x^n_2$ and $x_1,x_2$, we have now two more terms to estimate, which are the following ones:
\begin{gather}
\notag
\left[\mathcal{P}_\mathrm{\Lambda}\left(\alpha\sum_{i=1}^2 z_i\left(f_i\left(\abs{\cdot-x_i}\right)-f_i\left(\abs{\cdot-x^n_i}\right)\right)*\frac{1}{|\cdot|}\right)\mathcal{P}_\mathrm{\Lambda},P^0\right],\\
\notag
\left[\mathcal{P}_\mathrm{\Lambda}\left(\alpha\sum_{i=1}^2 z_i\left(f_i\left(\abs{\cdot-x_i}\right)-f_i\left(\abs{\cdot-x^n_i}\right)\right)*\frac{1}{|\cdot|}\right)\mathcal{P}_\mathrm{\Lambda},Q^n\right].
\end{gather}
Now, let us start noting that there exists a constant $\kappa>0$ such that for any $g\in\coulomb$
\[
\norm{g*\frac{1}{|\cdot|}}_{\schatten{\infty}}\leq\kappa E\left(\mathrm{\Lambda}\right)\norm{g}_\coulomb,
\]
where $E(\mathrm{\Lambda})=\sqrt{1+\mathrm{\Lambda}^2}$. This is an immediate consequence of the following inequality proved in \cite[Proof of Theorem~1, Step~3]{HaiLewSer-2005-JPA}
\[
\abs{g*\frac{1}{|\cdot|}}\leq\kappa\norm{g}_\coulomb \abs{D^0}.
\]
Thus,
\begin{multline*}
\norm{\left[\mathcal{P}_\mathrm{\Lambda}\left(\alpha\sum_{i=1}^2 z_i\left(f_i\left(\abs{\cdot-x_i}\right)-f_i\left(\abs{\cdot-x^n_i}\right)\right)*\frac{1}{|\cdot|}\right)\mathcal{P}_\mathrm{\Lambda},Q^n\right]}_{\schatten{2}}\\
\leq 2 E\left(\mathrm{\Lambda}\right)\alpha\kappa\left(z_1\norm{\overline{f}_1-\overline{f}^n_1}_\coulomb+z_2\norm{\overline{f}_2-\overline{f}^n_2}_\coulomb\right)\norm{Q^n}_{\schatten{2}}\\
\leq \underbrace{2 E\left(\mathrm{\Lambda}\right)\alpha\kappa}_{C_5}\left(z_1\norm{\overline{f}_1-\overline{f}^n_1}_\coulomb+z_2\norm{\overline{f}_2-\overline{f}^n_2}_\coulomb\right)\left(\norm{Q^n-Q}_{\schatten{2}}+\norm{Q}_{\schatten{2}}\right).
\end{multline*}
It is important to remark that in both the previous inequalities the constant $\kappa$ does not depend on $n$, and $\overline{f}_i$ and $\overline{f}^n_i$ are short notations for $f_i\left(\abs{\cdot-x_i}\right)$ and $f_i\left(\abs{\cdot-x^n_i}\right)$, $i=1,2$, respectively.\\
Moreover, let $\varphi^n=\mathcal{P}_\mathrm{\Lambda}\left(\alpha\sum_{i=1}^2 z_i\left(f_i\left(\abs{\cdot-x_i}\right)-f\left(\abs{\cdot-x^n_i}\right)\right)*\frac{1}{|\cdot|}\right)\mathcal{P}_\mathrm{\Lambda}$. Since the kernel of $\varphi^n P^0$ in the Fourier domain is $\left(2\pi\right)^{3/2}\widehat{\varphi^n}\left(p-q\right)P^0\left(q\right)$, we get
\begin{align*}
    \abs{\left[\widehat{\varphi^n},\widehat{P^0}\right]\left(p,q\right)}^2&=\left(2\pi\right)^3\abs{\widehat{\varphi^n}\left(p-q\right)}^2\tr_{\C^4}\left(P^0\left(p\right)-P^0\left(q\right)\right)^2\\
    &=2\left(2\pi\right)^3 \abs{\widehat{\varphi^n}\left(p-q\right)}^2 \tr_{\C^4}\left(P^0\left(p\right)P^0_\perp\left(q\right)\right),
\end{align*}
where $P^0_\perp=1-P^0$. By \cite[Lemma~12]{HaiLewSer-2005-CMP}, we have
\[
    \tr_{\C^4}\left(P^0\left(p\right)P^0_\perp\left(q\right)\right)\leq\frac{\abs{p-q}^2}{2E\left(\left(p+q\right)/2\right)^2},
\]
and thus
\[
    \norm{\left[\varphi_Q,P^0\right]}^2_{\schatten{2}}\leq \left(2\pi\right)^3\iint_{\abs{p},\abs{q}\leq\Lambda}\frac{\abs{\widehat{\alpha\left(z_1\left(\overline{f}_1-\overline{f}_1^n\right)+z_2\left(\overline{f}_2-\overline{f}_2^n\right)\right)}\left(p-q\right)}^2}{\abs{p-q}^2 E\left(\left(p+q\right)/2\right)^2}\diff p \diff q,
\]
which yields
\[
    \norm{\left[\varphi_Q,P^0\right]}_{\schatten{2}}\leq \alpha C_4 \left(z_1\norm{\overline{f}_1-\overline{f}_1^n}_\coulomb+z_2\norm{\overline{f}_2-\overline{f}_2^n}_\coulomb\right),
\]
where $C_4$ is some constant dependent on the cut-off parameter $\Lambda$.\\
Now, $\norm{\widehat{\rho}_{\left[\varphi^n,Q^n\right]}}_\coulomb=0$. Indeed,
\begin{multline*}
\widehat{\rho}_{\left[\varphi^n,Q^n\right]}\left(k\right)=\frac{1}{\left(2\pi\right)^{3/2}}\left[\int\int\widehat{\varphi^n}\left(p+\frac{k}{2}-s\right)\tr_{\C^4}\widehat{Q}^n\left(s,p-\frac{k}{2}\right)\diff p\diff s\right.\\
\left.-\int\int\widehat{\varphi^n}\left(s-p+\frac{k}{2}\right) \tr_{\C^4}\widehat{Q}^n\left(p+\frac{k}{2},s\right)\diff p\diff s \right]=0,
\end{multline*}
by the following variable transformation:
\begin{align*}
    \begin{cases}
    s=p'+\frac{k}{2},\\
    p=s'+\frac{k}{2}.
    \end{cases}
\end{align*}
Moreover $\rho_{\left[\varphi^n,P^0\right]}=0$ because the Dirac matrices $\alpha_i$, $i=1,2,3$, and $\beta$ have zero $\C^4$-trace.\\
Finally,
\[
\norm{Q^n-Q}\leq\int_0^t\left(\alpha_n+\beta_n\norm{Q^n-Q}\right)\diff s,
\]
where
\[
\alpha_n=\sup_{\substack{s\in\intcc{0}{\tau}\\\overline{x}_1,\overline{x}_2\in\R^3}} \left[\left(C_4+C_5\norm{Q\left(s\right)}\right)\left(z_1\norm{\overline{f}_1-\overline{f}^n_1}_\coulomb+z_2\norm{\overline{f}_2-\overline{f}^n_2}_\coulomb\right)\right]
\]
goes to zero as $n$ goes to infinity by \cref{stat:f-convergence} and
\begin{multline*}
    \beta_n=\sup_{\substack{s\in\intcc{0}{\tau}\\\overline{x}_1,\overline{x}_2\in\R^3}}\left[C_1+C_2\left(\norm{Q\left(s\right)}+\norm{Q^n\left(s\right)}\right)+C_3\right.\\
    \left.+C_5\left(z_1\norm{\overline{f}_1-\overline{f}^n_1}_\coulomb+z_2\norm{\overline{f}_2-\overline{f}^n_2}_\coulomb\right)\right]
\end{multline*}
is bounded. Thus, Gr\"{o}nwall's lemma allows us to conclude.
\endproof
\begin{rmk}
Notice that in the previous proof we get some rough estimates, since the upper bound is obtained by using the term $E\left(\mathrm{\Lambda}\right)$ which is linear in $\mathrm{\Lambda}$ whereas in QED we always expect a logarithmic dependence (and in particular divergence) on the cut-off parameter.
\end{rmk}

\subsection{Local existence for the Newton equation}
Let us now focus on the Cauchy problem of the remaining equations in \eqref{eq:2-nuclei},
\begin{equation}\label{eq:newton-subsystem}
	\begin{cases}
	m_k\frac{\diff^2}{\diff t^2}x_k\left(t\right)=-\nabla_{x_k}W_Q\left(t,x_1,x_2\right),\quad k=1,2\\
	x_k\left(0\right)=x^0_k,\quad\frac{\diff}{\diff t}x_k\left(0\right)=v^0_k,\quad k=1,2
	\end{cases}
\end{equation}
where we recall that
\[
W_Q\left(t,x_1,x_2\right)=\alpha \mathcal{D}\left(\rho_Q,\sum_{i=1}^2 z_i f_i\left(\abs{\cdot-x_i}\right)\right)-\alpha \int \int \frac{z_1 f_1\left(\abs{x-x_1}\right) z_2 f_2\left(\abs{y-x_2}\right)}{\abs{x-y}}\diff x\diff y.
\]
We immediately remark that we deal with a standard ODE and so we only have to pay attention to the source term $W_Q$. We are now ready to state and prove the following lemma:
\begin{lem}\label{stat:newton-existence}
Let $Q\in \mathcal{B}^e_\tau$. Then system~\eqref{eq:newton-subsystem} has a unique (local) solution in $\left(C^2\left(\intcc{0}{\tau},\R^3\right)\right)^2$ which is also in $\mathcal{B}^{n,1}_\tau\times\mathcal{B}^{n,2}_\tau$. Moreover, the map
\[
\mathcal{F}^\mathrm{NEW}:Q\in\mathcal{B}^e_\tau\longmapsto\left(x_1,x_2\right)\in \left(\mathcal{B}^{n,1}_\tau\cap C^2\left(\intcc{0}{\tau},\R^3\right)\right)\times\left(\mathcal{B}^{n,2}_\tau\cap C^2\left(\intcc{0}{\tau},\R^3\right)\right)
\]
is bounded and continuous.
\end{lem}
\proof
First of all, we remark that the nucleus-nucleus interaction term in $W_Q$ can be rewritten as a scalar product in the Coulomb space due to our assumptions. Secondly we immediately know that system~\eqref{eq:newton-subsystem} has a unique local solution by the Cauchy-Lipschitz theorem since a direct computation shows that $-\nabla_x W_Q$ is a well defined Lipschitz function.\\
Now we want to prove that the map $\mathcal{F}^\mathrm{NEW}:Q\in \mathcal{B}^e_\tau\mapsto\left(x_1,x_2\right)\in \left(\mathcal{B}^{n,1}_\tau\cap C^2\left(\intcc{0}{\tau},\R^3\right)\right)\times\left(\mathcal{B}^{n,2}_\tau\cap C^2\left(\intcc{0}{\tau},\R^3\right)\right)$ is bounded and continuous. First, the map $\mathcal{F}^\mathrm{NEW}$ is well defined. Indeed,
\begin{align*}
\norm{\frac{\diff x_k}{\diff t}}_{C^0\left(\intcc{0}{\tau},\R^3\right)}&\leq\abs{v^0_k}+\frac{\tau}{m_k}\sup_{\intcc{0}{\tau}\times\left(\R^3\right)^2}\abs{\nabla_{x_k}W_Q}\\
&\leq\abs{v^0_k}+\alpha \frac{\tau}{m_k}C_f\left(\norm{Q}_{C^0\left(\intcc{0}{\tau},\opspace\right)}+C_f\right)\\
&\leq\abs{v^0_k}+\alpha\frac{\tau}{m_k}C_f\left(\norm{Q_I}+C_f\right)\\
&\leq \abs{v^0_k}+1,\quad k=1,2
\end{align*}
where the last inequality is justified by \eqref{eq:inequality-2}.\\
Now let $\left(Q^n\right)_{n\in\N}\subseteq\mathcal{B}^e_\tau$ be a sequence converging to $Q\in\mathcal{B}^e_\tau$ in $\opspace$. We want to prove that
\[
\mathcal{F}^\mathrm{NEW}\left(Q^n\right)=\left(x^n_1,x^n_2\right)\xrightarrow{n\to\infty}\left(x_1,x_2\right)=\mathcal{F}^\mathrm{NEW}\left(Q\right).
\]
Let $\widetilde{x}^n_k=x^n_k-x_k$, $k=1,2$. Then, from
\begin{multline*}
    m_1 \frac{\diff^2}{\diff t^2}x^n_1\left(t\right)=\alpha \mathcal{D}\left(\rho_{Q^n\left(t\right)},z_1\nabla f_1\left(\abs{\cdot-x^n_1\left(t\right)}\right)\right)\\
    -\alpha \mathcal{D}\left(z_1\nabla f_1\left(\abs{\cdot-x^n_1\left(t\right)}\right),z_2 f_2\left(\abs{\cdot-x^n_2\left(t\right)}\right)\right)\\
    -\alpha \mathcal{D}\left(\rho_{Q\left(t\right)},z_1\nabla f_1\left(\abs{\cdot-x_1\left(t\right)}\right)\right)\\
    +\alpha \mathcal{D}\left(z_1\nabla f_1\left(\abs{\cdot-x_1\left(t\right)}\right),z_2 f_2\left(\abs{\cdot-x_2\left(t\right)}\right)\right),
\end{multline*}
we derive
\begin{align*}
    m_1 \frac{\diff^2}{\diff t^2}\widetilde{x}^n_1\left(t\right)&=\alpha\left[\mathcal{D}\left(\rho_{Q^n\left(t\right)-Q\left(t\right)},z_1 \nabla f_1\left(\abs{\cdot-\overline{x}^n_1\left(t\right)}\right)\right)\right.\\ \nonumber
    &\left.+\mathcal{D}\left(\rho_{Q\left(t\right)},z_1\left(\nabla f_1\left(\abs{\cdot-\overline{x}^n_1\left(t\right)}\right)-\nabla f_1\left(\abs{\cdot-\overline{x}_1\left(t\right)}\right)\right)\right)\right]\\ \nonumber
    &-\alpha\left[\mathcal{D}\left(z_1\left(\nabla f_1\left(\abs{\cdot-\overline{x}^n_1\left(t\right)}\right)-\nabla f_1\left(\abs{\cdot-\overline{x}_1\left(t\right)}\right)\right),z_2 f_2\left(\abs{\cdot-\overline{x}^n_2\left(t\right)}\right)\right)\right.\\ \nonumber
    &\left.+\mathcal{D}\left(z_1 \nabla f_1\left(\abs{\cdot-\overline{x}_1\left(t\right)}\right),z_2 \left(f_2\left(\abs{\cdot-\overline{x}^n_2\left(t\right)}\right)-f_2\left(\abs{\cdot-\overline{x}_2\left(t\right)}\right)\right)\right)\right],
\end{align*}
which implies
\[
\abs{m_1\frac{\diff^2}{\diff t^2}\widetilde{x}^n_1\left(t\right)}\leq \alpha A_n + \alpha B_n\left(\abs{\widetilde{x}^n_1\left(t\right)}+\abs{\widetilde{x}^n_2\left(t\right)}\right),
\]
where
\[
A_n=\sup_{t\in\intcc{0}{\tau}}\mathcal{D}\left(\rho_{Q^n\left(t\right)-Q\left(t\right)},z_1 \nabla f_1\left(\abs{\cdot-x^n_1\left(t\right)}\right)\right)
\]
and
\[
0\leq B_n\leq \norm{z_1 f_1\left(|\cdot|\right)}_{L^2\left(\R^3\right)}\left(\norm{\rho_Q}_\coulomb+2\norm{z_2 f_2\left(|\cdot|\right)}_\coulomb\right).
\]
Now $\left(B_n\right)_{n\in\N}$ is bounded since $\rho_Q\in\coulomb$ and $f_1,f_2\in L^2\left(\R^3\right)\cap\coulomb$ by hypotheses. Moreover, $\left(A_n\right)_{n\in\N}$ goes to zero when $n$ goes to infinity. Indeed, it is enough to prove that $Q\in\schatten{2}\mapsto \mathcal{D}\left(\rho_Q,\nabla f\left(|\cdot|\right)\right)\in\R^3$ is continuous. This is actually true since
\[
\abs{\mathcal{D}\left(\rho_Q,\nabla f\left(|\cdot|\right)\right)}\leq \norm{\nabla f\left(|\cdot|\right)}_\coulomb \norm{\rho_Q}_\coulomb
\]
and we know that $Q\in\schatten{2}\mapsto\rho_Q\in L^2\left(\R^3\right)\cap \coulomb
$ is continuous as explained for instance in \cite[(12)]{Sab-2014-AMRE}. Likewise, the same estimate holds for $\widetilde{x}^n_2$. Since $\widetilde{x}^n_k\left(0\right)=\frac{\diff\widetilde{x}^n_k}{\diff t}\left(0\right)=0$, we get that $\widetilde{x}^n_k$ goes to zero in $C^2\left(\intcc{0}{\tau},\R^3\right)$ when $n$ goes to infinity by Gr\"{o}nwall's lemma both for $k=1$ and for $k=2$.
\endproof

\subsection{Schauder fixed-point argument}
We proved so far local existence results for both \eqref{eq:BDF-subsystem} and \eqref{eq:newton-subsystem}. We are now ready to take advantage of our previous work, by applying a standard fixed-point argument in order to prove a local existence result for system~\eqref{eq:2-nuclei}.
\begin{thm}\label{stat:2-particles-existence}
System~\eqref{eq:2-nuclei} has a solution $\left(Q,\overline{x}_1,\overline{x}_2\right)$ in
\[
C^1\left(\intcc{0}{\tau},\opspace\right)\times\left(C^2\left(\intcc{0}{\tau},\R^3\right)\right)^2.
\]
\end{thm}
\proof
Let us denote by $i$ the compact injection
\[
i:\left(\mathcal{B}^{n,1}_\tau\cap C^2\left(\intcc{0}{\tau},\R^3\right)\right)\times\left(\mathcal{B}^{n,2}_\tau\cap C^2\left(\intcc{0}{\tau},\R^3\right)\right)\longrightarrow \mathcal{B}^{n,1}_\tau\times\mathcal{B}^{n,2}_\tau.
\]
We can now define the functional $\mathcal{K}=i\circ\mathcal{F}^\mathrm{NEW}\circ\mathcal{F}^\mathrm{BDF}$ which maps $\mathcal{B}^{n,1}_\tau\times\mathcal{B}^{n,2}_\tau$ into itself: if $\left(y_1,y_2\right)\in\mathcal{B}^{n,1}_\tau\times\mathcal{B}^{n,2}_\tau$, $z=\mathcal{K}\left(y_1,y_2\right)$ satisfies
\begin{equation}
	\begin{cases}
	i\frac{\diff}{\diff t}Q\left(t\right) = \left[D_{Q,y_1,y_2},Q\left(t\right)\right]+\left[V_{Q,y_1,y_2},P^0\right],\\[8pt]
	m\frac{\diff^2}{\diff t^2}z_k\left(t\right) = -\nabla_{x_k}W_Q\left(t,z_1,z_2\right),\quad k=1,2\\[8pt]
P\left(0\right) = P_I,\;P\left(t\right)^2 = P\left(t\right),\;Q\left(t\right) = P\left(t\right) - P\left(0\right) \in\opspace,\\[8pt]
	z_k\left(0\right) = x^0_k,\;\frac{\diff z_k}{\diff t}\left(0\right) = v^0_k,\quad k=1,2
\end{cases}
\end{equation}
with $\left(P,z_1,z_2\right)\in\mathcal{B}^e_\tau\times\mathcal{B}^{n,1}_\tau\times\mathcal{B}^{n,2}_\tau$. The map $\mathcal{K}$ is continuous and compact since $\mathcal{F}^\mathrm{NEW}$ and $\mathcal{F}^\mathrm{BDF}$ are continuous and bounded (by \cref{stat:newton-existence} and \cref{stat:BDF-existence} respectively) and the injection $i$ is continuous and compact. Since $\mathcal{B}^{n,1}_\tau\times\mathcal{B}^{n,2}_\tau$ is convex and bounded, by the Schauder fixed-point theorem $\mathcal{K}$ has a fixed point $\left(x_1,x_2\right)\in\mathcal{B}^{n,1}_\tau\times\mathcal{B}^{n,2}_\tau$ which is actually in $\left(\mathcal{B}^{n,1}_\tau\cap C^2\left(\intcc{0}{\tau},\R^3\right)\right)\times\left(\mathcal{B}^{n,1}_\tau\cap C^2\left(\intcc{0}{\tau},\R^3\right)\right)$ and $\left(P,x_1,x_2\right)$ is a solution of system~\eqref{eq:2-nuclei} with $P=\mathcal{F}^\mathrm{BDF}\left(x_1,x_2\right)$.
\endproof
We can now immediately state an important property of the solution of system~\eqref{eq:BDF-subsystem}, namely $P\left(t\right)=Q\left(t\right)+P^0$ is indeed a projection for any $t\in\intcc{0}{\tau}$.
\begin{lem}
Let $\left(Q,x_1,x_2\right)\in C^1\left(\intcc{0}{\tau},\opspace\right)\times \left(C^2\left(\intcc{0}{\tau},\R^3\right)\right)^2$ be given by \cref{stat:2-particles-existence}. Then $P\left(t\right)=Q\left(t\right)+P^0$ is an orthogonal projector.
\end{lem}
\proof
See \cite[Lemma~3.2]{HaiLewSpa-2005-LMP}.
\endproof
As a consequence, we deduce that $Q\left(t\right)=P\left(t\right)-P^0$ is $P^0$-trace class and its $P^0$-trace is conserved all along the time evolution. Indeed, by Lemma $2$ in \cite{HaiLewSer-2005-CMP} we know that $Q\left(t\right)=P\left(t\right)-P^0$ is $P^0$-trace class and that its charge $\tr_{P^0}\left(Q\left(t\right)\right)$ is an integer. Moreover,
\[
\tr_{P^0}\left(Q\left(t\right)\right)=\tr\left(Q\left(t\right)^3\right),\forall t\in\intco{0}{\tau}.
\]
Since $t\in\intco{0}{\tau}\mapsto Q\left(t\right)$ is continuous in the $\schatten{2}$ and therefore also in the $\schatten{3}$ topology, we get that $t\in\intco{0}{\tau}\mapsto \tr_{P^0}\left(Q\left(t\right)\right)$ is continuous and hence constant, namely
\[
\tr_{P^0}\left(Q\left(t\right)\right)=\tr_{P^0}\left(Q_I\right),\forall t\in\intco{0}{\tau}.
\]

\section{Uniqueness}\label{sec:uniqueness}
This section is devoted to prove a uniqueness result for the solution of system~\eqref{eq:2-nuclei}. Thus, we prove the following proposition, by a classical application of Gr\"{o}nwall's lemma:
\begin{prop}
The solution $\left(Q,x_1,x_2\right)$ to \eqref{eq:2-nuclei} is unique in
\[
C^1\left(\intoo{0}{\tau},\opspace\right)\times\left(C^2\left(\intcc{0}{\tau},\R^3\right)\right)^2.
\]
\end{prop}
\proof
Let $\left(Q,x_1,x_2\right)$ and $\left(Q',x'_1,x'_2\right)$ be two solutions of \eqref{eq:2-nuclei}. Let $\widetilde{x}_1=x_1-x'_1$, $\widetilde{x}_2=x_2-x'_2$ and $\widetilde{Q}=Q-Q'$. We define
\[
h\left(t\right)=\left(\abs{\widetilde{x}_1\left(t\right)}+\abs{\widetilde{x}_2\left(t\right)}+\norm{\widetilde{Q}\left(t\right)}_{\schatten{2}}\right)^p,
\]
where $p>2$. By means of a computation like in the proof of \cref{stat:newton-existence}, we get
\begin{align*}
    m_1 \frac{\diff^2}{\diff t^2}\widetilde{x}_1\left(t\right)&=\alpha\left[\mathcal{D}\left(\rho_{Q-Q'},z_1 \nabla f_1\left(\abs{\cdot-x_1}\right)\right)\right.\\
    &\left.+\mathcal{D}\left(\rho_{Q'},z_1\left(\nabla f_1\left(\abs{\cdot-x_1}\right)-\nabla f_1\left(\abs{\cdot-x'_1}\right)\right)\right)\right]\\
    &-\alpha\left[\mathcal{D}\left(z_1\left(\nabla f_1\left(\abs{\cdot-x_1}\right)-\nabla f_1\left(\abs{\cdot-x'_1}\right)\right),z_2 f_2\left(\abs{\cdot-x_2}\right)\right)\right.\\
    &\left.+\mathcal{D}\left(z_1 \nabla f_1\left(\abs{\cdot-x_1}\right),z_2 \left(f_2\left(\abs{\cdot-x_2}\right)-f_2\left(\abs{\cdot-x'_2}\right)\right)\right)\right]
\end{align*}
and an easy computation shows that
\[
\abs{\widetilde{x}_1\left(t\right)}\leq\int_0^t \left(t-s\right)\abs{\frac{\diff^2}{\diff s^2}\widetilde{x}_1\left(s\right)}\diff s\leq K_1\int_0^t\left(t-s\right)\left(\abs{\widetilde{x}_1\left(s\right)}+\abs{\widetilde{x}_2\left(s\right)}+\norm{\widetilde{Q}\left(s\right)}_{\schatten{2}}\right)\diff s,
\]
with $K_1$ constant dependent only on the cut-off parameter $\mathrm{\Lambda}$ and on the $L^2$ and Coulomb norms of $f_1,f_2$. The same estimate trivially holds for $\widetilde{x}_2\left(t\right)$. On the other hand,
\begin{align*}
\norm{\widetilde{Q}\left(t\right)}_{\schatten{2}}&\leq\int_0^t\norm{\frac{\diff}{\diff s}\left(Q-Q'\right)}_{\schatten{2}}\diff s=\int_0^t\norm{F\left(Q\right)-F\left(Q'\right)}_{\schatten{2}}\diff s\\
&\leq K_2 \int_0^t \left(\abs{\widetilde{x}_1\left(s\right)}+\abs{\widetilde{x}_2\left(s\right)}+\norm{\widetilde{Q}\left(s\right)}_{\schatten{2}}\right)\diff s,
\end{align*}
where here $f_i\left(|\cdot|\right)=f\left(\abs{\cdot-x_i}\right)$ and $f'_i\left(|\cdot|\right)=f\left(\abs{\cdot-x'_i}\right)$, $i=1,2$, and $K_2$ is again a constant dependent only on $\mathrm{\Lambda}$ and the norms of $f_1,f_2$ as $K_1$.\\
Now,
\begin{align*}
h\left(t\right)&\leq\left(\int_0^t\left(K_2+2K_1\left(t-s\right)\right)\left(\abs{\widetilde{x}_1\left(s\right)}+\abs{\widetilde{x}_2\left(s\right)}+\norm{\widetilde{Q}\left(s\right)}_{\schatten{2}}\right)\diff s\right)^p\\
&\leq \left(\int_0^t\left(K_2+2K_1\left(t-s\right)\right)^{p'}\diff s\right)^{\frac{p}{p'}}\left(\int_0^t h\left(s\right)\diff s\right),
\end{align*}
where $p'$ is the conjugate exponent associated to $p$, which implies that $h\left(t\right)=0$, $\forall t\in\intcc{0}{\tau}$.
\endproof
\begin{rmk}
    This work is deeply inspired by \cite{CanLeB-1999-MMMAS}, where nuclei and electrons are treated separately and then combined via the Schauder fixed-point theorem. For this reason, the proof of the local existence and uniqueness result of \cref{stat:main-result} follows the same argument. Alternatively, one could apply the Cauchy-Lipschitz theorem directly to the system \eqref{eq:M-nuclei} of coupled differential equations.
\end{rmk}

\section{Global-in-time existence}\label{sec:global-in-time-existence}
We recall that the energy of the system described by \eqref{eq:2-nuclei} is given by
\begin{multline*}
E^{\left(2\right)}\left(Q\left(t\right),x_1\left(t\right),x_2\left(t\right)\right)\\
=\tr_{P^0}\left(D^0Q\left(t\right)\right)-\alpha \mathcal{D}\left(\rho_{Q\left(t\right)},\sum_{k=1}^2 z_k f_k\left(\abs{\cdot-x_k\left(t\right)}\right)\right)+\frac{\alpha}{2}\mathcal{D}\left(\rho_{Q\left(t\right)},\rho_{Q\left(t\right)}\right)\\
-\frac{\alpha}{2}\int\int\frac{\abs{Q\left(x,y\right)}^2}{\abs{x-y}}\diff x\diff y+\frac{1}{2}\sum_{k=1}^2 m_k \abs{\dot{x}_k\left(t\right)}^2\\
+\alpha \int \int \frac{z_1 f_1\left(\abs{x-x_1\left(t\right)}\right) z_2 f_2\left(\abs{y-x_2\left(t\right)}\right)}{\abs{x-y}}\diff x\diff y.
\end{multline*}\\
It is easy to prove that $E^{\left(2\right)}$ is conserved along any solution of system~\eqref{eq:2-nuclei}. Indeed,
\begin{prop}\label{stat:energy-conservation}
Let $\left(Q,x_1,x_2\right)\in C^1\left(\intcc{0}{\tau},\opspace\right)\times \left(C^2\left(\intcc{0}{\tau},\R^3\right)\right)^2$ be given by \cref{stat:2-particles-existence}. Then $E^{\left(2\right)}\left(Q\left(t\right),x_1\left(t\right),x_2\left(t\right)\right)=E^{\left(2\right)}\left(Q_I,x^0_1,x^0_2\right)$ for any $t\in\intcc{0}{\tau}$.
\end{prop}
\proof
\begin{multline*}
\frac{\diff}{\diff t}E^{\left(2\right)}\left(Q\left(t\right),x_1\left(t\right),x_2\left(t\right)\right)\\
=\tr_{P^0}\left(D^0\dot{Q}\left(t\right)\right)+\alpha \mathcal{D}\left(\rho_{Q\left(t\right)}-\sum_{k=1}^2 z_k f_k\left(\abs{\cdot-x_k\left(t\right)}\right),\rho_{\dot{Q}\left(t\right)}\right)-\alpha\tr\left(\frac{Q\left(x,y\right)}{\abs{x-y}}\dot{Q}\left(t\right)\right)\\
+\sum_{k=1}^2 \scalprod{\dot{x}_k\left(t\right)}{m_k \Ddot{x}_k\left(t\right)+\nabla_{x_k}W_Q\left(t,x_1,x_2\right)}.
\end{multline*}
We consequently obtain
\[
\frac{\diff}{\diff t}E^{\left(2\right)}\left(Q\left(t\right),x_1\left(t\right),x_2\left(t\right)\right)=\tr_{P^0}\left(D_{Q\left(t\right),x_1\left(t\right),x_2\left(t\right)}\dot{Q}\left(t\right)\right).
\]
If we now insert \cref{eq:BDF-subsystem}, we get
\begin{multline*}
    \frac{\diff}{\diff t}E^{\left(2\right)}\left(Q\left(t\right),x_1\left(t\right),x_2\left(t\right)\right)=-i\tr_{P^0}\left(D_{Q\left(t\right),x_1\left(t\right),x_2\left(t\right)}\left[D_{Q\left(t\right),x_1\left(t\right),x_2\left(t\right)},Q\left(t\right)\right]\right)\\
    -i\tr_{P^0}\left(D_{Q\left(t\right),x_1\left(t\right),x_2\left(t\right)}\left[V_{Q\left(t\right),x_1\left(t\right),x_2\left(t\right)},P^0\right]\right)
\end{multline*}
and a simple computation shows that the two terms on the right-hand side are equal to zero. Thus the time derivative of the energy vanishes and our statement is proved.
\endproof
Thanks to the conservation of energy of the system, we are now able to prove a global-in-time existence result for the unique solution of \eqref{eq:2-nuclei}.
\begin{cor}
Let $\left(Q,x_1,x_2\right)\in C^1\left(\intcc{0}{\tau_\mathrm{max}},\opspace\right)\times\left(C^2\left(\intcc{0}{\tau_\mathrm{max}},\R^3\right)\right)^2$ be the maximal solution provided by \cref{stat:2-particles-existence}. If $0\leq\alpha<\frac{4}{\pi}$, then $\tau_\mathrm{max} =+\infty$.
\end{cor}
\proof
By \cite{BacBarHelSie-1999-CMP}, we know the following inequality:
\begin{multline*}
E^{\left(2\right)}\left(Q\left(t\right),x_1\left(t\right),x_2\left(t\right)\right)+\frac{\alpha}{2}\sum_{k=1}^2 \norm{z_k f_k\left(\abs{\cdot-x_k\left(t\right)}\right)}_\coulomb^2\\
\geq\left(1-\alpha\frac{\pi}{4}\right)\tr_{P^0}\left(D^0Q\left(t\right)\right)+\frac{\alpha}{2}\norm{\rho_{Q\left(t\right)}-\sum_{k=1}^2 z_k f_k\left(\abs{\cdot-\overline{x}_k\left(t\right)}\right)}_\coulomb^2\\
+\frac{1}{2}\sum_{k=1}^2 m_k\abs{\dot{x}_k\left(t\right)}^2.
\end{multline*}
By \cref{stat:energy-conservation}, $E^{\left(2\right)}\left(Q\left(t\right),x_1\left(t\right),x_2\left(t\right)\right)=E^{\left(2\right)}\left(Q_I,x^0_1,x^0_2\right)$. Moreover, both $\abs{\dot{x}_k}$ and $\tr_{P^0}\left(D^0Q\left(t\right)\right)$ are bounded since $\norm{f\left(\abs{\cdot-x_k\left(t\right)}\right)}_\coulomb=\norm{f\left(|\cdot|\right)}_\coulomb$, $k=1,2$.\\
Finally we get that $\abs{x_k\left(t\right)}$ is bounded too and $\norm{Q\left(t\right)}_{\schatten{2}}$ stays bounded for $\alpha<\pi/4$ since $\tr_{P^0}\left(D^0Q\right)\geq\norm{Q}_{\schatten{2}}$. Therefore $\tau_\mathrm{max}=+\infty$.
\endproof

\addtocontents{toc}{\protect\setcounter{tocdepth}{0}} 

\section*{Acknowledgements}
The author would like to thank \'Eric Séré and William Borrelli for useful remarks and fruitful discussions. This research project has received funding from the European Union's Horizon 2020 research and innovation programme under the Marie Skłodowska-Curie grant agreement N°945332.

\section*{Conflict of interest}
The author declares no conflict of interest.

\section*{Data availability statement}
Data sharing is not applicable to this article as it has no associated data. 

\addtocontents{toc}{\protect\setcounter{tocdepth}{2}} 

\bibliographystyle{siam} 
\small{\bibliography{bibliography}}

\end{document}